\documentclass[11pt]{{article}}
\usepackage{CJK}
\usepackage{amsmath}
\usepackage{amsfonts}
\usepackage{mathrsfs}
\usepackage{amssymb,amsmath}
\usepackage{amstext}
\usepackage{array}
\usepackage{theorem}
\usepackage{latexsym}
\usepackage[all]{xy}
\usepackage{amscd}
\usepackage{bbding}
\usepackage{graphics}
\usepackage{graphicx}
\usepackage{bbm}
\usepackage{color}
\usepackage{pgf,pgfarrows,pgfnodes,pgfautomata,pgfheaps}
  
\usepackage[all]{xy}
\def\Tr{\mathop{\rm Tr}\nolimits}
\def\mod{\mathop{\rm mod}\nolimits}

\title{\Large \bf A Note On Gorenstein Projective Conjecture II
}
\author{ Xiaojin Zhang\thanks{{\it E-mail address}:
xjzhang@nuist.edu.cn}\\
{\footnotesize \it School of Mathematics and Statistics, Nanjing
University of Information Science and Technology},\\{\footnotesize
\it Nanjing 210044, P. R. China}}
\date{}

\begin{document}
 \begin{CJK*}{GBK}{song}

\baselineskip=18pt
 \maketitle
 \noindent{\bf Abstract}\ \ In this paper, we prove that
Gorenstein projective conjecture is left and right symmetric and the
co-homology vanishing condition can not be reduced in general.
Moreover, Gorenstein projective conjecture is proved to be true for
CM-finite algebras.


\noindent{\bf AMS(2000) Subject Classification}\ \  16G10, 16E05.
\vspace{0.5cm}

\noindent{\large \bf 1\ \ Introduction}

\vspace{0.2cm}

For the representation theory of Artinian algebras, the
Auslander-Reiten conjecture (ARC) which is related to generalized
Nakayama conjecture (GNC) is everything. It was proposed by
Auslander and Reiten, which says that $M$ is projective if ${\rm
Ext}_{\Lambda}^{i}(M\bigoplus\Lambda,M\bigoplus\Lambda)=0$ for any
$i\geq1$(See [2,3]). Achievements for special cases have been got by
K. R. Fuller, B. Zimmermann-Huisgen, A. Mar$\dot{o}$ti and G.
Wilson...(See [10,15,16]). In general it is still open now. As a
special case of Auslander-Reiten conjecture, Luo and Huang proposed
the following Gorenstein projective conjecture (GPC) in 2008:

 Let $\Lambda$ be an Artinian algebra and let $M$ be a Gorenstein
projective module. Then $M$ is projective if and only if ${\rm
Ext}_{\Lambda}^{i}(M,M)=0$ for any $i\geq1$(See [14]).

It is still unknown whether the Auslander-Reiten conjecture is left
and right symmetric. But as we stated Gorenstein projective
conjecture is a special case. So what about the left and right
symmetric property of Gorenstein projective conjecture? In this
paper, we will give a positive answer to this question.

By the definition of Gorenstein projective conjecture, for an
algebra $\Lambda$ the truth of Auslander-Reiten conjecture implies
the truth of Gorenstein projective conjecture. So we can get a large
class of algebras satisfying Gorenstein projective conjecture. It is
interesting to ask: Is there an algebra satisfying Gorenstein
projective conjecture while for which the Auslander-Reiten
conjecture is unknown?

Recall that an algebra is called CM-finite (of finite Cohen-Macaulay
type) if there are only finitely many isomorphism classes of
indecomposable finitely generated Gorenstein projective modules.
CM-finite algebras are studied by several authors recently (see
[5,6,7,12,13]). Although the Auslander-Reiten conjecture for this
class of algebras is unknown, we will give a positive answer to the
second question above.

The paper is organized as follows:

In Section 2, based on some facts of Gorenstein projective modules,
we will show the symmetric property of Gorenstein projective
conjecture. Moreover, an example is given to show that the condition
'${\rm Ext}_{\Lambda}^{i}(M,M)=0$ for any $i\geq1$' in Gorenstein
projective conjecture can not be reduced to '${\rm
Ext}_{\Lambda}^{i}(M,M)=0$ for some positive integer $t$ and any
$1\leq i\leq t$ .'

In Section 3, we will prove that CM-finite algebras satisfy the
Gorenstein projective conjecture by showing the Gorenstein
projective conjecture holds for algebras with finite self-orthogonal
indecomposable Gorenstein projective modules (up to isomorphisms).

Throughout the paper, $\Lambda$ is an Artinian algebra and all
modules are finitely generated left $\Lambda$-modules.

\vspace{0.5cm}

\noindent{\large \bf 2\ \ Symmetric property of Gorenstein
projective conjecture }

\vspace{0.2cm}

In this section we will show the symmetric property of Gorenstein
projective conjecture. First, we need to recall some notions and
lemmas. The following definition is due to Auslander, Briger, Enochs
and Jenda (see [1,8,9]).

{\bf Definition 2.1} A module $M$ is called {\it Gorenstein
projective} if for any $i\geq1$

(1) ${\rm Ext}_{\Lambda}^{i}(M,\Lambda)=0$

(2) ${\rm Ext}_{\Lambda}^{i}(\Tr M,\Lambda)=0$

Where $\Tr M$ denotes the Auslander transpose of $M$.

 Let $\dots\rightarrow P_2(M)\rightarrow P_1(M)\rightarrow M\rightarrow 0$ be a minimal projective resolution of $M$.
 Denoted by $\Omega^i M$ the $i$-th syzygy of $M$. Dually, one can define $\Omega^{-i}M$. We remark that for any $i\geq 0$ $\Omega^i M$ is
a Gorenstein projective if so is $M$. Let $\mathscr{C}$ be the
subcategory of $\mod\Lambda$ consisting of modules $M$ such that
${\rm Ext}_{\Lambda}^{j}(M,\Lambda)=0$ for any $j\geq 1$ and
$\mathscr{D}$ a subcategory of $\mod\Lambda$ consisting of
Gorenstein projective modules. We use $\underline{\mathscr{C}}$ and
$\underline{\mathscr{D}}$ to denote the stable subcategory of
$\mathscr{C}$ and $\mathscr{D}$ modulo projectives, respectively. We
recall the following proposition from [1].

\vspace{0.2cm}

{\bf Proposition 2.2} (1) {\it
$\Omega:\underline{\mathscr{C}}\rightarrow \underline{\mathscr{C}}$
is a fully faithful functor .}

(2) {\it $\Omega:\underline{\mathscr{D}}\rightarrow
\underline{\mathscr{D}}$ is an equivalence.}

(3) {\it $(-)^*={\rm Hom}(-,\Lambda):\mathscr{D}\rightarrow
\mathscr{D}^o$ is a duality, where $\mathscr{D}^o$ denotes the
subcategory of Gorenstein projective right $\Lambda$-modules.}

{\it Proof.} (1) is a result of Auslander and Bridger. One can get
(2) by (1) and the remark above. (3) is well-known.
$\hfill{\square}$

\vspace{0.2cm} Recall that a module $M$ is called self-orthogonal if
${\rm Ext}_{\Lambda}^{j}(M, M)=0$ for any $j\geq 1$. The following
self-orthogonal property is essential to the main result in this
section.

{\bf Lemma 2.3} {\it Let $\Lambda$ be an algebra. Then for any
$M\in\mathscr{D}$ and $i\geq1$, $M$ is self-orthogonal if and only
if $M^*$ is self-orthogonal}.

\vspace{0.2cm}

{\it Proof.} Since $(-)^*$ is a duality between $\mathscr{D}$ and
$\mathscr{D}^o$, it is enough to show that ${\rm
Ext}_{\Lambda}^{i}(M,M)=0$ implies ${\rm
Ext}_{\Lambda}^{i}(M^*,M^*)=0$.

  One can take the following minimal projective resolution of $M$:
$$\dots\rightarrow P_1\rightarrow P_0\rightarrow M\rightarrow0
\ \ \ \ (1)$$ Applying the functor ${\rm Hom}(-,M)$ to sequence
$(1)$ above, since ${\rm Ext}_{\Lambda}^{i}(M,M)=0$ we get the
following exact sequence $$0\rightarrow {\rm
Hom}(M,M)\rightarrow{\rm Hom}(P_0,M)\rightarrow{\rm
Hom}(P_1,M)\rightarrow\cdots\ \ \ (2)$$

On the other hand, applying the functor $(-)^*$ to the sequence
$(1)$, since $M\in\mathscr{D}\subseteq\mathscr{C}$ one can show the
following exact sequence $$0\rightarrow M^*\rightarrow
{P_0}^*\rightarrow{P_1}^*\rightarrow\cdots\ \ \ (3)$$ Then by using
the functor ${\rm Hom}(M^*,-)$ on the exact sequences (3), one has
the following exact sequence $$0\rightarrow {\rm
Hom}(M^*,M^*)\rightarrow{\rm Hom}(M^*,{P_0}^*)\rightarrow{\rm
Hom}(M^*,{P_1}^*)\rightarrow\cdots\ \ \ (4)$$ Using Proposition
2.2(3), we get ${\rm Ext}_{\Lambda}^{i}(M^*,M^*)\simeq {\rm
Ext}_{\Lambda}^{i}(M,M)=0$ by comparing sequences (2) with (4).
$\hfill{\square}$

\vspace{0.2cm}

Although the symmetric property of Auslander-Reiten conjecture is
still unknown, we are able to show the symmetric properties of
Gorenstein projective conjecture.

{\bf Theorem 2.4} {\it Let $\Lambda$ be an algebra and let
$\Lambda^o$ be the opposite ring of $\Lambda$. Then $\Lambda$
satisfies the Gorenstein projective conjecture if and only if
$\Lambda^o$ satisfies the Gorenstein projective conjecture.}

{\it Proof.} $\Rightarrow$ Assume that $N\in\mathscr{D}^o$ and ${\rm
Ext}_{\Lambda}^{i}(N,N)=0$ for any $i\geq 1$. By Proposition 2.2,
there is a $M\in\mathscr{D}$ such that $M^*\simeq N$. By Lemma 2.3
one gets that ${\rm Ext}_{\Lambda}^{i}(M,M)=0$. Note that $\Lambda$
satisfies the Gorenstein projective conjecture, we have $M$ is
projective, and hence $N\simeq M^*$ is projective.  Conversely, one
can formula the proof above. $\hfill{\square}$

\vspace{0.2cm}

Notice that Gorenstein projective conjecture is a special case of
Auslander-Reiten conjecture. It is natural to consider whether the
assumption of Gorenstein projective conjecture can be reduced. In
particular, whether can the condition '${\rm
Ext}_{\Lambda}^{i}(M,M)=0$ for any $i\geq1$' in GPC be reduced to
'${\rm Ext}_{\Lambda}^{i}(M,M)=0$ for some positive integer $t$ and
any $1\leq i\leq t$ ? At the end of this section, we construct an
example to give a negative answer to the question.
 \vspace{0.2cm}

{\bf Example 2.5} Let $n>t+1$ be a positive integer and let
$\Lambda$ be the algebra generated by the following quiver
 $$
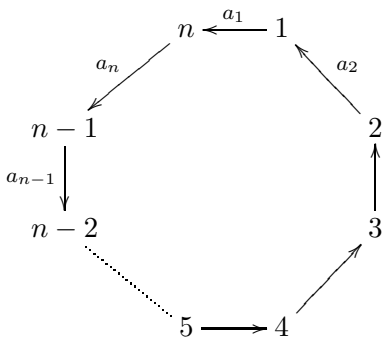
\xymatrix{
&n\ar[ld]_{a_n}&1\ar[l]_{a_1}&\\
n-1\ar[d]_{a_{n-1}}&&&2\ar[ul]_{a_2}\\
n-2&&&3\ar[u]\\
&5\ar@{.}[lu]\ar[r]&4\ar[ur]& }
$$ modulo the ideal $\{a_{n}a_{1}=0,a_{i}a_{i+1}=0|1\leq i\leq
n-1\}$. Denoted by $S(j)$ the simple module  according to the dot
$j$. Then

 (1) $\Lambda$ is a Nakayama self-injective algebra.

 (2) $S(j)$ is Gorenstein projective such that ${\rm
Ext}_{\Lambda}^{i}(S(j),S(j))=0$ for $t\geq i\geq1$ and $1\leq j
\leq n$, but it is not projective.

\vspace{0.5cm}

\noindent{\large \bf 3\ \ Gorenstein projective conjecture for
CM-finite algebras }

\vspace{0.2cm}

In this section we try to find a class of algebras which satisfy
Gorenstein projective conjecture and for which the Auslander-Reiten
conjecture is unknown. They are CM-finite algebras. We begin with
the following definition due to Beligiannis

{\bf Definition 3.1} An algebra is called CM-finite (of finite
Cohen-Macaulay type ) if there are only finite number of
indecomposable Gorenstein projective modules (up to isomorphisms).

\vspace{0.2cm}

{\bf Remark 3.2} (1) Algebras of finite representation type or
finite global dimension are CM-finite.

(2) There does exist a CM-finite algebra $\Lambda$ such that
$\Lambda$ is of infinite type and the global dimension of $\Lambda$
is infinite [13].

(3) There does exist a CM-finite algebra which is not Gorenstein
[5].

(4)An algebra with a trivial maximal n-orthogonal subcategory for
some positive integer $n$ is CM-finite [11].

\vspace{0.2cm}

Let $ \mathscr{C}$, $ \mathscr{D}$ and $\mathscr{D}^o$ be as in
Section 2. The following lemma partly from [1] plays an important
role in the proof of the main results.

{\bf Lemma 3.3} {\it For any $M\in \mathscr{C}$ and
$N\in\mod\Lambda$, then ${\rm Ext}_{\Lambda}^{1}(M,N)\simeq
\underline{{\rm Hom}}_{\Lambda}(\Omega^1 M, N)$ and hence ${\rm
Ext}_{\Lambda}^{i}(M,N)\simeq \underline{{\rm
Hom}}_{\Lambda}(\Omega^i M, N)$ for any $i\geq1$.}

{\it Proof.} The first assertion is a result of Auslander and
Bridger. For the second one, the case $i=1$ is clear. We only need
to show the case $i\geq 2$. Taking a minimal projective resolution
of $M$, one gets ${\rm Ext}_{\Lambda}^{i}(M,N)\simeq {\rm
Ext}_{\Lambda}^{1}(\Omega^{i-1}M,N) $ for any $i\geq2$. Notice that
$M\in\mathscr{C}$, by Proposition 2.2 one can show
$\Omega^{i-1}M\in\mathscr{C}$. Using the first assertion, we are
done. $\hfill{\square}$.

\vspace{0.2cm}

 The following
proposition gives a connection between the self-orthogonal property
of $M$ and that of $\Omega^i M$ for any $i\geq 0$.

{\bf Proposition 3.4 } {\it Let $M\in \mathscr{C}$ $(\mathscr{D})$.
Then

(1) $\Omega^i M$ is self-orthogonal in $\mathscr{C}$ $(\mathscr{D})$
for any $i\geq 0$ if $M$ is self-orthogonal.

(2) If $M\in\mathscr{D}$ is self-orthogonal, then ${\rm Tr}M$ is
self-orthogonal in $\mathscr{D}^o.$}

{\it Proof.} (1) For the case $i=0$, there is nothing to prove. By
Proposition 2.2, we only need to prove the case $i=1$, that is,
${\rm Ext}_{\Lambda}^{j}(\Omega M, \Omega M)=0$ for any $j\geq1$.
One gets ${\rm Ext}_{\Lambda}^{j}(\Omega M, \Omega M)\simeq
\underline{{\rm Hom}}_{\Lambda}(\Omega^{j+1} M, \Omega M)\simeq
\underline{{\rm Hom}}_{\Lambda}(\Omega^{j} M, M)$ by Proposition 2.2
and Lemma 3.3. Using the second equation of Lemma 3.3, one can show
$\underline{{\rm Hom}}_{\Lambda}(\Omega^{j} M, M)\simeq {\rm
Ext}_{\Lambda}^{j}(M,M)=0$ since M is self-orthogonal.

(2) Taking a minimal projective resolution of $M$, it is not
difficult to show that ${\rm Tr}M\simeq (\Omega^2 M)^*$ since
$M\in\mathscr{D}$. By Propositions 2.2 and 3.4(1), $\Omega^2 M$ is
also self-orthogonal in $\mathscr{D}$. Then by Lemma 2.3 and
Proposition 2.2 one gets the assertion. $\hfill{\square}$

\vspace{0.2cm}

The following proposition is crucial to the main results.

{\bf Proposition 3.5} {\it Let $\Lambda$ be an algebra with only
finite (up to isomorphism) self-orthogonal indecomposable modules in
$\mathscr{D}$ ($\mathscr{C}$) and let $M$ be a self-orthogonal
indecomposable module in $\mathscr{D}$ ($\mathscr{C}$). Then $M$ is
projective.}

{\it Proof.} Denoted by $\{M_1,M_2,...,M_t\}$ the complete set of
non-isomorphic self-orthogonal indecomposable modules in
$\mathscr{D}$ ($\mathscr{C}$). Then $M\simeq M_i$ for some $1\leq
i\leq t$.

Suppose that M is not projective. Then by Proposition 2.2 and
Proposition 3.4, we have the following set of self-orthogonal
indecomposable modules $\mathcal {S}=\{\Omega^i M|1\leq i \}$ in
$\mathscr{D}$ ($\mathscr{C}$).

We claim that there are two modules $\Omega^i M,\Omega^j M$ in
$\mathcal {S}$ such that $\Omega^i M\simeq\Omega^j M$ for some $i<
j$. Otherwise, one gets infinite number of non-isomorphic
self-orthogonal indecomposable modules in $\mathscr{D}$
($\mathscr{C}$), a contradiction. Again by Proposition 2.2, one gets
$M\simeq \Omega^{j-i}M$.

 Considering the following exact sequence $0\rightarrow
\Omega^{j-i}M\rightarrow P\rightarrow \Omega^{j-i-1}M\rightarrow 0
$, we will show ${\rm
Ext}_{\Lambda}^{1}(\Omega^{j-i-1}M,\Omega^{j-i}M)=0$. Since
$\Omega^{j-i-1}M\in \mathscr{D}$ ($\mathscr{C}$) and $M\simeq
\Omega^{j-i}M$, we get ${\rm
Ext}_{\Lambda}^{1}(\Omega^{j-i-1}M,\Omega^{j-i}M)\simeq
\underline{{\rm Hom}}_{\Lambda}(\Omega^{j-i} M, M)\simeq {\rm
Ext}_{\Lambda}^{j-i}(M,M)=0 $ by Proposition 2.2 and Lemma 3.3, and
hence $M$ is projective, a contradiction. The assertion holds
true.$\hfill\square$

\vspace{0.2cm}

 Now we are in the position to show the main result of
this section.

{\bf Theorem 3.6} {\it Let $\Lambda$ be CM-finite. Then $\Lambda$
satisfies Gorenstein projective conjecture.}

{\it Proof. } Since $\Lambda$ is CM-finite, then there are only
finite (up to isomorphisms) indecomposable modules in $\mathscr{D}$.
One can show the result by Proposition 3.5. $\hfill{\square}$

Although the Auslander-Reiten conjecture for CM-finite algebras is
unknown now, we have the following

\vspace{0.2cm}

{\bf Proposition 3.7} {\it Let $\Lambda$ be a CM-finite algebra and
let $M$ be a $\Lambda$-module satisfying ${\rm
Ext}_{\Lambda}^{i}(M,M\bigoplus\Lambda)=0$ for any $i\geq 1$. Then
the following are equivalent.

(1) $M$ is projective.

(2) $M$ is Gorenstein projective.}

{\it Proof.} $(1)\Rightarrow(2)$ is trivial. The converse follows
from Theorem 3.6. $\hfill{\square}$

\vspace{0.2cm}

We end this section with two open questions related to this paper.

{\bf Question 1} Does the Gorenstein projective conjecture hold for
virtually Gorenstein algebras (see [5])?

{\bf Question 2} Does the Auslander-Reiten conjecture hold for
CM-finite algebras?

\vspace{0.2cm}

{\bf Acknowledgements}  Part of the paper was finished when the
author stayed in University of Bielefeld with the support of DAAD
Fellowship. The author would like to  thank people in Bielefeld for
their help. He also wants to thank Prof. Xiaowu Chen, Prof. Zhaoyong
Huang, Prof. Shengyong Pan and Prof. Changchang Xi for useful
suggestions. The research of the author is supported by NSFC(Grant
No.11101217) and NSF for Colleges and Universities in Jiangsu
Province (Grant No.11KJB110007).

 \end{CJK*}
\end{document}